\documentclass[12pt,twoside]{article}
\usepackage{graphicx}
\usepackage{amsfonts}
\usepackage{amsbsy}
\usepackage{amssymb}
\usepackage{amsmath,tikz,pgfplots}
\usetikzlibrary{decorations.pathmorphing}
\usepackage{cite}
\usepackage{graphics}
\def\bpsp{\begin{pspicture}}
\def\epsp{\end{pspicture}}

\newtheorem{theorem}{Theorem}[section]
\newtheorem{remark}[theorem]{Remark}
\newtheorem{example}[theorem]{Example}
\newtheorem{lemma}[theorem]{Lemma}
\newtheorem{corollary}[theorem]{Corollary}
\newtheorem{definition}[theorem]{Definition}
\newtheorem{proposition}[theorem]{Proposition}
\newtheorem{note}{Note}
\newtheorem{case}{Case}
\newtheorem{conjecture}{Conjecture}
\newtheorem{question}{Question}

\newcommand{\bea}{\begin{eqnarray}}
\newcommand{\eea}{\end{eqnarray}}
\newcommand{\beq}{\begin{eqnarray*}}
\newcommand{\eeq}{\end{eqnarray*}}

\def\m4{\mbox{\rm ~(mod $4$)}}

\def \bd{\begin{definition}}
\def \ed{\end{definition}}
\def \bqu{\begin{question}}
\def \equ{\end{question}}
\def \bcc{\begin{conjecture}}
\def \ecc{\end{conjecture}}
\def \bt{\begin{theorem}}
\def \et{\end{theorem}}
\def \bl{\begin{lemma}}
\def \el{\end{lemma}}
\def \bc{\begin{corollary}}
\def \ec{\end{corollary}}
\def \be{\begin{equation}}
\def \ee{\end{equation}}
\def \ben{\begin{enumerate}}
\def \een{\end{enumerate}}
\def \ba{\begin{array}}
\def \ea{\end{array}}
\def \bp{\begin{proposition}}
\def \ep{\end{proposition}}
\def \bx{\begin{example}}
\def \ex{\end{example}}
\def \br{\begin{remark}}
\def \er{\end{remark}}
\def \bdsc{\begin{description}}
\def \edsc{\end{description}}

\def \bn{\begin{case}}
\def \en{\end{case}}
\def \bnt{\begin{note}}
\def \ent{\end{note}}
\def\1{1\!\!1}

\def\mm2{\mbox{\rm ~(mod $2$)}}
\def\m4{\mbox{\rm ~(mod $4$)}}

\def\qed{\nolinebreak\hfill\rule{.2cm}{.2cm}\par\addvspace{.5cm}}

\def\m{\mu}

\def\1{\textbf{1}}
\def\0{\textbf{0}}

\begin{document}
\title{On the Estrada index of unicyclic and bicyclic signed graphs }
\author{ Tahir Shamsher$^{a}$, S. Pirzada$^{b}$, Mushtaq A. Bhat$^{c}$ \\
$^{}$ {\em  $^{a,b}$Department of Mathematics, University of Kashmir, Srinagar, Kashmir, India}\\
$^{}${\em  $^{c}$Department of Mathematics, National Institute of Technology, Srinagar, India}\\
$ ^{a} $\texttt{tahir.maths.uok@gmail.com}; $^{b}$\texttt{pirzadasd@kashmiruniversity.ac.in}\\
 $ ^{a} $\texttt{mushtaqab@nitsri.net}}
\date{}

\pagestyle{myheadings} \markboth{ Tahir, Pirzada, Bhat}{On the Estrada index of unicyclic and bicyclic signed graphs }
\maketitle
\vskip 5mm
\noindent{\footnotesize \bf Abstract. } Let  $\Gamma=(G, \sigma)$ be a signed graph of order $n$ with eigenvalues $\mu_1,\mu_2,\ldots,\mu_n.$  We define  the Estrada index of a signed graph $\Gamma$ as $EE(\Gamma)=\sum_{i=1}^ne^{\mu_i}$.  We characterize the signed unicyclic graphs with the maximum Estrada index. The signed graph $\Gamma$ is said to have the pairing property if $\mu$ is an eigenvalue whenever $-\mu$ is an eigenvalue of $\Gamma$ and both $\mu$ and $-\mu$ have the same multiplicities. If $\Gamma_{p}^-(n, m)$ denotes the set of all unbalanced graphs on $n$ vertices and $m$ edges with the pairing property, we determine the signed graphs having the maximum Estrada index in $\Gamma_{p}^-(n, m)$, when $m=n$ and $m=n+1$. Finally, we find the signed graphs among all unbalanced complete bipartite signed graphs having the maximum Estrada index.

\vskip 3mm

\noindent{\footnotesize Keywords: Signed graph, Estrada index, pairing property, largest eigenvalue, spectral moment.  }

\vskip 3mm
\noindent {\footnotesize AMS subject classification:  05C22, 05C50.}

\section{Introduction}\label{sec1}

A signed graph $\Gamma$ is an ordered  pair $(G, \sigma)$ in which $G$ is an underlying graph and $\sigma$ is a function from the edge set $E(G$) to $\{-1,1\}$, which is called
a signing. For a signed graph $\Gamma=(G, \sigma)$ and its subgraph $H \subset G$, we use the notation   $(H, \sigma)$ for  writing the signed subgraph of $\Gamma=(G, \sigma)$, where $\sigma$ is the restriction of the mapping $\sigma: E(G) \rightarrow \{-1,1\}$ to the edge set $E(H)$. The adjacency matrix $A_{\Gamma}=\left(a_{i j}\right)$ of a signed graph $\Gamma=(G, \sigma)$ naturally arose from the unsigned graph by putting $-1$ or $1$, whenever the corresponding edge is either negative or positive, respectively. The characteristic polynomial, denoted by $\varphi(\Gamma, x)=\operatorname{det}\left(x I-A_{\Gamma}\right)$, is called the characteristic polynomial of the signed graph $\Gamma=(G, \sigma)$. For brevity, the spectrum of the adjacency matrix $A_{\Gamma}$ is called the spectrum of the signed graph $(G, \sigma)$. Let the signed graph $\Gamma$ of order $n$ has distinct eigenvalues $\mu_1(\Gamma),\mu_2(\Gamma),\ldots,\mu_k(\Gamma)$ (we drop $\Gamma$ where the signed graph is understood) and let their respective multiplicities be $m_1,m_2,\ldots,m_k$. The adjacency spectrum of $\Gamma$ is written as $Spec(\Gamma)=\{\mu^{(m_1)}_1(\Gamma),\mu^{(m_2)}_2(\Gamma),\ldots,\mu^{(m_k)}_k(\Gamma)\}$. From the definition, it follows that the matrix $A_{\Gamma}$ is a real symmetric and hence the eigenvalues $\mu_{1}((G, \sigma)) \geq \mu_{2}((G, \sigma)) \geq \dots \geq \mu_{n}((G, \sigma))$ of the signed graph $(G, \sigma)$ are all real numbers. The largest eigenvalue $\mu_{1}(\Gamma)$ is also known as  the index of the signed graph $\Gamma$.\\
\indent The concept of signature switching is necessary when dealing with signed graphs. Let $Z$ be a  subset of the vertex set $V(\Gamma)$. The switched signed graph $\Gamma^{Z}$ is obtained from $\Gamma$ by reversing the signs of the edges in the cut $[Z, V(\Gamma) \backslash Z]$. Clearly, we see that the signed graphs $\Gamma$ and $\Gamma^{Z}$ are switching equivalent. The switching equivalence is an equivalence relation that preserves the eigenvalues. The switching class of $\Gamma$ is denoted by $[\Gamma]$. The sign of a cycle is the product of the signs of its edges. A signed cycle $C^{}_{n\sigma}$ on $n$ vertices is  positive (or negative) if it contains an even (or odd) number of negative edges, respectively. A signed graph is said to be balanced if all of its signed cycles are positive, otherwise, it is unbalanced. That is, a signed graph is said to be balanced if it switches to the signed graph with all positive signature. Otherwise, it is said to be unbalanced.  By $\sigma \sim+,$ we say that the signature $\sigma$ is equivalent to the all-positive signature, and the corresponding signed graph is equivalent to its underlying graph. In general, the signature is determined by the set of positive cycles. Hence, all trees are switching equivalent to the all-positive signature. Equivalently, we can say that the edge signs of bridges are irrelevant. Finally, the signature of the balanced signed graphs is switching equivalent to the all-positive one. In our drawings of signed graphs, we represent the negative edges with dashed lines and the positive edges with solid lines. A connected signed graph is said to be unicyclic if it has the same number of vertices and edges. If the number of edges is one more than the number of vertices, then it is said to be bicyclic. For definitions and notations of graphs, we refer to \cite{sp}.\\
\indent The signed graph $\Gamma$ is said to have the pairing property if $\mu$ is an eigenvalue whenever $-\mu$ is an eigenvalue of $\Gamma$ and both $\mu$ and $-\mu$ have the same multiplicities. The signed graph $\Gamma=(G,+)$ with all positive signature has the pairing property if and only if its underlying graph $G$ is bipartite. For any signature $\sigma$ it is not true. \\
\indent The rest of the paper is organized as follows. In Section $2$, we extend the Estrada index to signed graphs. Furthermore, we characterize the signed unicyclic graphs with the maximum Estrada index. In Section $3$,  we find  the signed graphs in the set of all unbalanced  unicyclic and bicyclic graphs having the pairing property  with the maximal Estrada index respectively.

\section{Estrada index in signed graphs}\label{sec2}

The Estrada index, a graph-spectrum-based structural descriptor, of a graph is defined as the trace of the adjacency matrix exponential and was first proposed by Estrada in $2000$. Pena et al. \cite{de} recommended calling it the Estrada index, which has since become widely adopted. This index can be used to measure a range of things, including the degree of protein folding \cite{e1, e2, e3}, the subgraph centrality and bipartivity of complex networks \cite{e4,e5}.  Because of the graph Estrada index's exceptional use,  various Estrada indices based on the eigenvalues of other graph matrices were investigated. Estrada index-based invariant concerning the Laplacian matrix, signless Laplacian matrix, distance matrix,   distance Laplacian matrix and distance signless Laplacian matrix have been studied, see \cite{e6}.\\
In social networks, the balance (stability) \cite{e7,e8} of a signed network
can be quantified by
\begin{equation}
k=\frac{tr(e^{A_{(G,\sigma)}})}{tr(e^{A_{(G,+)}})},
\end{equation}
 where $tr(X)$ denotes the trace of the matrix $X$. Motivated by Equation $(2.1)$, we define the Estrada index for the signed graph  in full analogy with the Estrada index for a  graph, as
\begin{equation}
EE(\Gamma)=EE((G,\sigma))=\sum_{i=1}^ne^{\mu_i},
\end{equation}
where $\mu_1$, $\mu_2$, $\ldots$, $\mu_n$ are the eigenvalues of the signed graph $\Gamma$. The Seidel matrix $S_G$ of a simple graph $G$ with $n$ vertices and having the adjacency matrix $A_G$ is defined as
$S_G = J -I -2A_G$. Obviously, the Seidel matrix is the adjacency matrix of some signed graph $\Gamma=(K_n, \sigma)$, where $K_n$ is the complete graph on $n$ vertices. Therefore, Eq. $(2.2)$ is the extension of the Siedel Estrada index \cite{se}. \\
For non-negative integer $k$, let $M_{k}(\Gamma)=\sum_{i=1}^{n} \mu_{i}^{k}$  denote the $k$-th spectral moment of $\Gamma$.  From the Taylor expansion of $\mathrm{e}^{\mu_{i}}$,  $E E(\Gamma)$ in $(2.2)$ can be rewritten as
\begin{equation}
E E(\Gamma)=\sum_{k=0}^{\infty} \frac{M_{k}(\Gamma)}{k !}.
\end{equation}
It is well-known that $M_{k}(\Gamma)$ is equal to the difference of the number of positive and negative closed walks of length $k$ in $\Gamma$. Mathematically, we have
\begin{equation}
M_{k}(\Gamma)=w^+(k)-w^-(k),
\end{equation}
where $w^+(k)$ and $w^-(k)$ are, respectively, the number of positive and negative closed walks of length $k$ in $\Gamma$. In particular, we have
\begin{center}

$ M_{0}(\Gamma)=n$, $~ M_{1}(\Gamma)=0$,  $~ M_{2}(\Gamma)=2m$ and $~ M_{3}(\Gamma)=6(t^+-t^-)$,
\end{center}
where  $n$ is the number of vertices,  $m$ is the number of edges, $t^+$ is the number of positive triangles and $t^-$ is the number of negative triangles in the signed graph $\Gamma$.

Let $\Gamma_{1}$ and $\Gamma_{2}$ be two signed graphs. If $M_{ k}\left(\Gamma_{1}\right) \geq M_{ k}\left(\Gamma_{2}\right)$ holds for any positive integer $k$, then, by Eq. $(2.3)$   we get $E E\left(\Gamma_{1}\right) \geq E E\left(\Gamma_{2}\right)$. Further, if the strict inequality $M_{k}\left(\Gamma_{1}\right)>M_{ k}\left(\Gamma_{2}\right)$ holds for at least one integer $k$, then $E E\left(\Gamma_{1}\right)>E E\left(\Gamma_{2}\right)$.
It is easy to see that if $\Gamma$ has $q$ connected components $\Gamma_{1}, \Gamma_{2}, \ldots, \Gamma_{q}$, then $E E(\Gamma)=$ $\sum_{i=1}^{q} E E\left(\Gamma_{i}\right).$ So, we shall investigate the Estrada index of connected signed graph from now on. One classical problem of graph spectra is to identify the extremal graphs with respect to the Estrada index in some given class of graphs, for example, see \cite{h,zu, g1}. For a signed tree, all signatures are equivalent. The following result shows that among all signed trees on $n$ vertices, the signed path $P_n$ has a minimum and the signed star $S_n$ has the maximum Estrada index.
\begin{figure}
\centering
	\includegraphics[scale=.8]{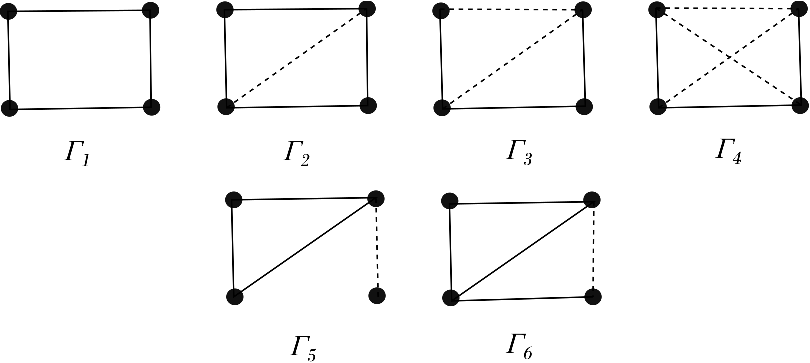}
	\caption{Signed graph $\Gamma_i \quad i=1,2,\ldots,6.$}
	\label{Figure 1}
\end{figure}

\begin{theorem} {\em \cite{h}}
 If $T_{n}$ is an $n$-vertex tree different from $S_{n}$ and $P_{n}$, then
$$
EE\left((P_{n}, \sigma\right))<EE\left((T_{n}, \sigma\right))<EE\left((S_{n}, \sigma\right)).
$$
\end{theorem}
 {\bf Remark 1.1} Let $\Gamma = (G, +)$ be a connected signed graph of order $n$ with all positive signature and $e$ a positive edge. The signed graph $\Gamma^{\prime}=\Gamma+e$ is obtained from $\Gamma$ by adding the edge $e$. It is easy to see that any self-returning walk of length $k$ of $\Gamma$ is also a self-returning walk of length $k$ of $\Gamma^{\prime}$. Thus, $M_{ k}\left(\Gamma'\right) \geq M_{ k}\left(\Gamma\right)$ and
$EE(\Gamma^{\prime})\geq EE(\Gamma).$
 But in general adding any signed edge between two non-adjacent vertices of the signed graph $\Gamma = (G, \sigma)$ may not increase the Estrada index. Consider the signed graphs $\Gamma_i$, $i=1, ~2, \ldots, 6$ as shown in Figure $1$. Their spectrum is given by $Spec(\Gamma_1)=\{2,-2, 0^2\}$, $Spec(\Gamma_2)=\{\frac{-1+\sqrt{17}}{2},1, 0, \frac{-1-\sqrt{17}}{2}\}$, $Spec(\Gamma_3)=\{2,1,-1,-2\}$, $Spec(\Gamma_4)=\{\sqrt{5},1,-1, -\sqrt{5}\}$,  $Spec(\Gamma_5)=\{2.17,0.31,-1,-1.48\}$ and $Spec(\Gamma_6)=\{2,1,-1,-2\}$ respectively. The signed graph $\Gamma_2$ is obtained from $\Gamma_1$ by adding  negative edge between two non-adjacent vertices and clearly $EE(\Gamma_2)=8.55<9.52= EE(\Gamma_1).$ The signed graph $\Gamma_4$ is obtained from $\Gamma_3$ by adding  negative edge between two non-adjacent vertices and  $EE(\Gamma_4)=12.54 > 10.61= EE(\Gamma_3).$ Furthermore, The signed graph $\Gamma_6$ is obtained from $\Gamma_5$ by adding  positive edge  and  $EE(\Gamma_5)=10.71 > 10.61= EE(\Gamma_6).$ Therefore, edge addition (deletion) technique cannot be used to compare the Estrada index in signed graphs.
\begin{theorem}\label{2.3}
Let $G$ be a graph on $n$ vertices. Then $EE((G,+))\geq EE((G,-))$, with strict inequality if and only if $G$ contains at least one odd cycle.
\end{theorem}
{\bf Proof.}  Let $G$ be a graph on $n$ vertices. Put  $\Gamma_1=(G,+)$ and $\Gamma_2=(G,-)$. Then, by Eq. $(2.2)$, we have
\begin{equation}
\begin{aligned}
EE(\Gamma_1)-EE(\Gamma_2)&=\sum_{i=1}^ne^{\mu_i(\Gamma_1)}-\sum_{i=1}^ne^{\mu_i(\Gamma_2)} \\
&= \sum_{i=1}^n(e^{\mu_i(\Gamma_1)}-e^{\mu_i(\Gamma_2)}).
\end{aligned}
\end{equation}
The signed graph $\Gamma_2$ can be obtained from the signed graph $\Gamma_1$ by negating each edge. Therefore $Spec(\Gamma_1)=-Spec(\Gamma_2)$. Thus, by rearrangement of Eq. $(2.5)$ and using Taylor's expansion, we have
\begin{equation}
\begin{aligned}
EE(\Gamma_1)-EE(\Gamma_2)&= \sum_{i=1}^n(e^{\mu_i(\Gamma_1)}-e^{-\mu_i(\Gamma_1)})\\
&= 2\sum_{k=0}^{\infty}\frac{M_{2k+1}(\Gamma_1)}{(2k+1)!}.
\end{aligned}
\end{equation}
The signed graph $\Gamma_1$ has all positive signature and therefore by Eq. $(2.4)$, $M_{2k+1}(\Gamma_1)\geq 0$. If $\Gamma_1$ has an odd cycle of size $l$, then $M_{l}(\Gamma_1)> 0$.  Hence the result follows. \qed
There exist exactly two switching classes on the signings of an odd unicyclic graph (whose unique cycle has odd girth). The cycle with all positive  signature and the cycle with all negative signature. The following result is an immediate consequence of  Theorem \ref{2.3}.

\begin{corollary} \label{1.4}
Let $G$ be  an odd unicyclic graph of order $n$ and let $\Gamma_1$ be any balanced signed graph on $G$ and $\Gamma_2$ be any unbalanced one. Then $EE(\Gamma_1) > EE(\Gamma_2)$.
\end{corollary}

Let $G$ be a bipartite unicyclic graph of girth $l$ and order $n$ and let $\Gamma_1$ be any balanced signed graph on $G$ and $\Gamma_2$ be any unbalanced one. It is easy to see that $M_{2k+1}(\Gamma_1)=0=M_{2k+1}(\Gamma_2)$ for each $k\geq 0$,  $M_{2k}(\Gamma_1)=M_{2k}(\Gamma_2)$ for  $2k\leq l-2$ and $M_{2k}(\Gamma_1)>M_{2k}(\Gamma_2)$ for $2k\geq l$. In particular, $M_{l}(\Gamma_1)=M_{l}(\Gamma_2)+4l.$  Thus, by Eq.s $(2.3)$ and $(2.4)$, we have the following lemma.

\begin{lemma} \label{1.5}
Let $G$ be a bipartite unicyclic graph of order $n$ and let $\Gamma_1$ be any balanced signed graph on $G$ and $\Gamma_2$ be any unbalanced one. Then $EE(\Gamma_1) > EE(\Gamma_2)$.
\end{lemma}

\indent Let $\Gamma^+(n,l)$ and $\Gamma^-(n,l)$ denote the set of balanced and unbalanced unicyclic  graphs with $n$ vertices and containing a cycle of length $l\le n$ respectively. Also, we denote by $\Gamma_n^{l+}$ (respt. $\Gamma_n^{l-}$),  the signed graph obtained by identifying the center of the signed star $S_{n-l+1}$ with a vertex of positive cycle $C_{l+}$ (respt. negative cycle $C_{l-}$). Du et al. \cite{zu} characterized the unique signed unicyclic graph having all positive signature with the maximum Estrada index and showed the following.

\begin{lemma}\label{1.6}
 Let $\
 \Gamma=(G,+)$ be a unicyclic graph on $n \geq 4 $ vertices. Then $EE(\Gamma) \leq EE(\Gamma_n^{3+})$ with equality
if and only if $\Gamma$ is isomorphic to $\Gamma_n^{3+}$.
\end{lemma}

The following result is directly obtained from Corollary \ref{1.4}, Lemma \ref{1.5} and Lemma \ref{1.6}.

\begin{theorem}
 Let $\
 \Gamma=(G,\sigma)$ be a signed unicyclic graph on $n \geq 4 $ vertices. Then  $EE(\Gamma) \leq EE(\Gamma_n^{3+})$ with equality
if and only if $\Gamma$ is switching equivalent to $\Gamma_n^{3+}$.
\end{theorem}

Next, we show that the Estrada index of the balanced cycle $C_{n+}$ is almost equal to the Estrada index of the unbalanced cycle ${ C}_{n-}$.

\begin{theorem}
Let $C_{n+}$ and $C_{n-}$ be the balanced and unbalanced cycles on $n$ vertices respectively. Then $E E\left(C_{n-}\right) \approx n J_{0}\approx E E\left(C_{n+}\right),$
where $J_{0}=\sum_{r \geqslant 0} \frac{1}{(r !)^{2}}=2.27958530 \ldots . $ Also, $EE\left(C_{n+}\right)=EE\left(C_{n-}\right)+\epsilon_n$, where $\epsilon_n\rightarrow 0$ as $n\rightarrow \infty$.
\end{theorem}
{\bf Proof.} The Estrada index of the $n$-vertex signed cycle $C_{n+}$  can be approximated \cite{ec} as
\begin{equation}
E E\left(C_{n+}\right) \approx n J_{0},
\end{equation}
where $J_{0}=\sum_{r \geqslant 0} \frac{1}{(r !)^{2}}=2.27958530 \ldots . $\\
The eigenvalues of the unbalanced cycle  $C_{n-}$ are given by $2 \cos \frac{ (2r+1) \pi}{n}, \quad r=0, 1,2, \ldots, n-1$.
 Therefore
$$E E\left(C_{n-}\right)=\sum_{r=0}^{n-1} \mathrm{e}^{2 \cos ( (2r+1) \pi / n)} .$$
 The angles $ (2r+1) \pi / n$, for $r=0,1,2, \ldots, n-1,$  uniformly cover the semi-closed interval $[0,2 \pi)$. Now, using the property of definite  integrals as a sum, we have
$$
E E\left(C_{n-}\right)=n\left(\frac{1}{n} \sum_{r=0}^{n-1} \mathrm{e}^{2 \cos ((2 r+1) \pi / n)}\right) \approx n\left(\frac{1}{2 \pi} \int_{0}^{2 \pi} \mathrm{e}^{2 \cos x} \mathrm{~d} x\right).
$$
As $\mathrm{e}^{2 \cos x}$ is an even function, therefore
$$
\int_{0}^{2 \pi} \mathrm{e}^{2 \cos x} \mathrm{~d} x=2 \int_{0}^{\pi} \mathrm{e}^{2 \cos x} \mathrm{~d} x=\pi J_{0},
$$
where $J_{0}$  is a special value of the function encountered in the theory of Bessel function and can be seen in \cite{bes} as
$$
J_{0}=\sum_{r=0}^{\infty} \frac{1}{(r !)^{2}}=2.27958530 \ldots.
$$
In view of this,
\begin{equation}
E E\left(C_{n-}\right) \approx n J_{0}=2.27958530 n .
\end{equation}
Eqs.  $(2.7)$ and $(2.8)$ give
\begin{equation}
EE\left(C_{n-}\right) \approx n J_{0} \approx EE\left(C_{n+}\right) .
\end{equation}
Note that  $M_{k}(C_{n+})=M_{k}(C_{n-})$ for  $k\leq n-1$ and $M_{k}(C_{n+})\geq M_{k}(C_{n-})$ for $k\geq n$. In particular, $M_{n}(C_{n+})=M_{n}(C_{n-})+4n.$  Thus, by Eq.s $(2.3)$ and $(2.4)$, we get
\begin{equation} EE(C_{n+})-EE(C_{n-})= \frac{4n}{n!}+\sum_{k=n+1}^{\infty}\frac{M_{k}(C_{n+})-M_{k}(C_{n-})}{k!}.
\end{equation}
The eigenvalues of  $C_{n+}$ and $C_{n-}$ are, respectively, given by $2 \cos \frac{ 2(r-1) \pi}{n}$ and $2 \cos \frac{ (2r+1) \pi}{n}, \quad r=0, 1,2, \ldots, n-1.$ Therefore
\begin{equation}
\begin{aligned}
\sum_{k=n+1}^{\infty}\frac{M_{k}(C_{n+})-M_{k}(C_{n-})}{k!}&=\sum_{k=n+1}^{\infty}\frac{\sum_{r=0}^{n-1}2^k(\cos^k \frac{ 2(r-1) \pi}{n}- \cos^k \frac{ (2r+1) \pi}{n})}{k!}.
\end{aligned}
\end{equation}
The maximum value of the function $\cos^k \frac{ 2(r-1) \pi}{n}- \cos^k \frac{ (2r+1) \pi}{n}$ is $2$. Thus, by Eq. $(2.11)$, we have
\begin{equation}
\begin{aligned}
\sum_{k=n+1}^{\infty}\frac{M_{k}(C_{n+})-M_{k}(C_{n-})}{k!} &\leq \sum_{k=n+1}^{\infty}\frac{\sum_{r=0}^{n-1}2^{k+1}}{k!} \\ &=\sum_{k=n+1}^{\infty}\frac{2^{k+1}n}{k!} \\ &= \frac{2^{n+2}n}{(n+1)!}+ \frac{2^{n+3}n}{(n+2)!}+\frac{2^{n+4}n}{(n+3)!}+\cdots \\&\leq  \frac{2^{n+2}n}{n!}\left(\frac{1}{(n+1)}+ \frac{2}{(n+1)^2}+\frac{2^2}{(n+1)^3}+\cdots \right).
\end{aligned}
\end{equation}
The series $\frac{1}{(n+1)}+ \frac{2}{(n+1)^2}+\frac{2^2}{(n+1)^3}+\dots$ is an infinite geometric progression with common ratio $\frac{2}{(n+1)}$. By inequality $(2.12)$, we obtain
\begin{equation}
\begin{aligned}
\sum_{k=n+1}^{\infty}\frac{M_{k}(C_{n+})-M_{k}(C_{n-})}{k!}  &\leq  \frac{2^{n+2}n}{n!(n-1)} .
\end{aligned}
\end{equation}
Eq.s $(2.10)$ and $(2.13)$ imply that
\begin{equation}
\begin{aligned}
EE(C_{n+})-EE(C_{n-})  &\leq   \frac{2^{n+2}n+4n(n-1)}{n!(n-1)} ,
\end{aligned}
\end{equation}
where the term $ \frac{2^{n+2}n+4n(n-1)}{n!(n-1)}\sim \frac{2^{n+2}}{n!} $ tends to zero as the girth $n$ becomes large enough. Moreover, the accuracy  of the Eq. $(2.9)$  can be seen from the data given in Table 1. As seen from this data, except for the first few values of $n$ $(n\leq 9)$, the accuracy is more
than sufficient.\\
{\bf Table 1.} Approximate and exact values of the Estrada index of the $n$-vertex signed cycles $\left(C_{n+}\right)$ and $\left(C_{n-}\right).$
\begin{center}
\begin{tabular}{rrrr}
\hline $n$ & $E E\left(C_{n+}\right)$ & $n J_{0}$ & $E E\left(C_{n-}\right)$  \\
\hline
3 & $8.1248150$ & $6.8387561$ & $5.571899$ \\
4 & $9.5243914$ & $9.1183414$ & $8.7127342$  \\
5 & $11.4961863$ & $11.3979268$ & $11.2993665$ \\
6 & $13.6967139$ & $13.6775122$ & $13.658309$\\
7 & $15.9602421$ & $15.9570975$ & $15.9533523$  \\
8 & $18.2371256$ & $18.2366829$ & $18.2368574$ \\
9 & $20.5163225$ & $20.5162683$ & $20.5163962$  \\
10 & $22.7958591$ & $22.7958536$ & $22.7958491$  \\
11 & $25.0754389$ & $25.0754390$ & $25.0754200$  \\
12 & $27.3550237$ & $27.3550243$ & $27.3550195$  \\
13 & $29.6346089$ & $29.6346097$ & $29.6345864$  \\
14 & $31.9141942$ & $31.9141951$ & $31.9141892$  \\
15 & $34.1937795$ & $34.1937804$ & $34.1937780$ \\
\hline
\end{tabular}
\end{center}
Hence the result follows. \qed
The main tool used to prove Lemma \ref{1.6} is the construction of  mappings which increases the $k$-th spectral moment for each $k$ and using Eq.$(2.3)$. For example consider  the following result.

 \begin{theorem} {\em \cite{zu}} For all positive signature $\sigma$ and $4\leq l \leq n$, we have $EE(\Gamma_n^{(l+1)+})<EE(\Gamma_n^{l+})$.
\end{theorem}
But in signed unicyclic graphs, in general, we cannot construct the  mapping which increases the $k$-th spectral moment for each $k$. To defend this statement, consider the signed unicyclic graphs $\Gamma_5^{4-}$ and $\Gamma_5^{5-}$. Their spectra is given by $Spec(\Gamma_5^{5-})=\{\frac{1+\sqrt{5}}{2}^2, \frac{1-\sqrt{5}}{2}^2,-2\}$  and $Spec(\Gamma_5^{4-})=\{\sqrt{3},\sqrt{2}, 0,-\sqrt{2},-\sqrt{3}\}$, respectively. It is easy to see that not only $EE(\Gamma_5^{5-})=11.30> 11.18 = EE(\Gamma_5^{4-})$ but also $M_4(\Gamma_5^{5-})= 30 > 26= M_4(\Gamma_5^{4-})$ and $M_5(\Gamma_5^{5-})= -10 < 0= M_5(\Gamma_5^{4-})$. Thus the problem of finding the unbalanced unicyclic graphs with extremal Estrada index is interesting.

\section{Unbalanced  unicyclic and bicyclic graphs with the pairing property and with maximal Estrada index }\label{sec2}

In this section, we first characterize the unbalanced bipartite unicyclic graphs with the maximal Estrada index.  Given two non-increasing real number sequences $\alpha=\{\alpha_1,\alpha_2, \alpha_3,\ldots, \alpha_n\} $ and $\beta=\{\beta_{1}, \beta_{2}, \beta_{3}, \ldots, \beta_{n}\} $, we say that $\alpha$ majorizes $\beta$, denoted by $\alpha \succeq \beta$, if
$
 \sum_{j=1}^{t} \alpha_{j} \geq \sum_{j=1}^{t} \beta_{j}
$
for each $t=1, \ldots, n$, with equality for $t=n.$ Also, if $\alpha \neq \beta$ then $\alpha \succ \beta$.

\begin{lemma}\label{3.1}{\em \cite{sf}}
Let $f:\mathbb{R}\rightarrow \mathbb{R}$  be a strictly convex function. If $\alpha \succeq \beta$, then $\sum_{i=1}^{n} f\left(\alpha_{i}\right) \geq$ $\sum_{i=1}^{n} f\left(\beta_{i}\right) .$ Also, if $\alpha \neq \beta$, then $\sum_{i=1}^{n} f\left(\alpha_{i}\right)>\sum_{i=1}^{n} f\left(\beta_{i}\right)$.
\end{lemma}

Let $\Gamma_{p}(n, m)$ denote the set of all signed graphs on $n$ vertices and $m$ edges with the pairing property. The following result will be useful in the sequel.

\begin{theorem}\label{3.2}
Let $\Gamma_{1}$, $\Gamma_{2}$ $\in \Gamma_{p}(n, m)$  be two signed graphs on $n$ vertices and $m$ edges with the pairing property. If $\Gamma_{1}$ has exactly four non-zero eigenvalues and $\Gamma_{2}$ has at least four non-zero eigenvalues with $\mu_{1}\left(\Gamma_{1}\right)>$ $\mu_{1}\left(\Gamma_{2}\right)$, then $E E\left(\Gamma_{1}\right)>E E\left(\Gamma_{2}\right)$.
\end{theorem}
{\bf Proof.} Let $\Gamma_{1}$, $\Gamma_{2}$ $\in \Gamma_{p}(n, m)$. Therefore, we have $\sum_{i=1}^{n} \mu_{i}^{2}(\Gamma_1)=\sum_{i=1}^{n} \mu_{i}^{2}(\Gamma_2)=2 m$. The signed graphs $\Gamma_{1}$ and $\Gamma_{2}$ have the pairing property, so we get $M_{2k+1}(\Gamma_1)=0=M_{2k+1}(\Gamma_2)$ for each $k\geq 0$. Let $\mu_{1}(\Gamma_2)\geq \mu_{2}(\Gamma_2)\geq\mu_{3}(\Gamma_2)\geq \dots \geq \mu_{2r}(\Gamma_2)$, $r\geq 2$, be the non-zero eigenvalues of $\Gamma_2$. Thus, by Eq. $(2.3)$ and using the pairing property, we have
\begin{equation}
EE\left(\Gamma_{1}\right)= n-2r+ 2\sum_{k=0}^{\infty} \frac{g_{k}\left(\mu_{1}^{2}\left(\Gamma_{1}\right), \mu_{2}^{2}\left(\Gamma_{1}\right),0,\dots,0\right)}{(2 k) !}
\end{equation}
and
\begin{equation}
EE\left(\Gamma_{2}\right)= n-2r+ 2\sum_{k=0}^{\infty} \frac{g_{k}\left(\mu_{1}^{2}\left(\Gamma_{2}\right), \mu_{2}^{2}\left(\Gamma_{2}\right), \dots, \mu_{r}^{2}\left(\Gamma_{2}\right)\right)}{(2 k) !} ,
\end{equation}
where $g_k(x_1,x_2,x_3,\ldots,x_r)=x_1^k+x_2^k+x_3^k+\dots + x_r^k$, $k$ is a positive integer and $\mu_{1}^{2}\left(\Gamma_{1}\right)+\mu_{2}^{2}\left(\Gamma_{1}\right)=m=\mu_{1}^{2}\left(\Gamma_{2}\right)+\mu_{2}^{2}\left(\Gamma_{2}\right)+\cdots+\mu_{r}^{2}\left(\Gamma_{2}\right)$. Now, Eq.s $(3.15)$ and $(3.16)$ imply that
\begin{equation}
EE\left(\Gamma_{1}\right)-EE\left(\Gamma_{2}\right)= 2\sum_{k=2}^{\infty} \frac{[g_{k}\left(\mu_{1}^{2}\left(\Gamma_{1}\right), \mu_{2}^{2}\left(\Gamma_{1}\right),0,\ldots,0\right)-g_{k}\left(\mu_{1}^{2}\left(\Gamma_{2}\right), \mu_{2}^{2}\left(\Gamma_{2}\right), \dots, \mu_{r}^{2}\left(\Gamma_{2}\right)\right)]}{(2 k) !} .
\end{equation}
  We know that the function $f(x)=x^{2k}$ is strictly convex for every positive integer $k.$\\
As $\mu_{1}\left(\Gamma_{1}\right)>\mu_{1}\left(\Gamma_{2}\right)$ and $\sum_{i=1}^{n} \mu_{i}^{2}(\Gamma_1)=\sum_{i=1}^{n} \mu_{i}^{2}(\Gamma_2)=2 m$, therefore  the vector\\ $\alpha=(\mu_{1}^{2}\left(\Gamma_{1}\right), \mu_{2}^{2}\left(\Gamma_{1}\right),0,\ldots,0)$ majorizes $\beta =(\mu_{1}^{2}\left(\Gamma_{2}\right), \mu_{2}^{2}\left(\Gamma_{2}\right), \cdots, \mu_{r}^{2}\left(\Gamma_{2}\right)),$ that is $\alpha \succ \beta$. Thus, by Lemma \ref{3.1}, we have
$$
g_{k}\left(\mu_{1}^{2}\left(\Gamma_{1}\right), \mu_{2}^{2}\left(\Gamma_{1}\right),0,\ldots,0\right)>g_{k}\left(\mu_{1}^{2}\left(\Gamma_{2}\right), \mu_{2}^{2}\left(\Gamma_{2}\right), \dots, \mu_{r}^{2}\left(\Gamma_{2}\right)\right),
$$
for $k \geq 2 .$ Hence, by Eq. $(3.17)$, we get $E E\left(\Gamma_{1}\right)>E E\left(\Gamma_{2}\right)$ and the result follows.  \qed

 \begin{lemma}\label{3.3}{\em \cite{uc}}
Let $\Gamma^-(n,l)$ be the set of  unbalanced unicyclic  graphs with $n$ vertices and containing a cycle of length $l\le n$. Then \\
$(i)$ for any $\Gamma \in \Gamma^-(n,l)$, we have $\mu_1(\Gamma_n^{l-})\geq \mu_1(\Gamma)$ with equality if and only if $\Gamma$ is switching equivalent to $\Gamma_n^{l-}$;\\
$(ii)$ $\mu_1(\Gamma_n^{l-})>\mu_1(\Gamma_n^{(l+1)-}).$
\end{lemma}

\indent The following result \cite{horn} shows that the eigenvalues of an induced signed subgraph of the signed graph $\Gamma$ interlaces the eigenvalues of   $\Gamma$.

\begin{lemma}\label{3.4}
 Let $\Gamma$ be an $n$-vertex signed graph with eigenvalues $\mu_{1} \geq \dots \geq \mu_{n}$, and let $\Gamma^{\prime}$ be an induced signed subgraph of $\Gamma$ with $m$ vertices. If the eigenvalues of $\Gamma^{\prime}$ are $\lambda_{1} \geq \dots \geq \lambda_{m}$, then $\mu_{n-m+i} \leq \lambda_{i} \leq \mu_{i},$ where $i=1, \ldots, m .$
\end{lemma}
\indent We now recall Schwenk's formula and can be seen in \cite{b1}.
\begin{lemma}\label{3.5}
Let $u$ be any fixed vertex of a signed graph $\Gamma$. Then
$$
\varphi(\Gamma, x)=x \varphi(\Gamma-u, x)-\sum_{vu  \in E(\Gamma)} \varphi(\Gamma-v-u, x)-2 \sum_{Y \in \mathcal{Y}_{u}} \sigma(Y) \varphi(\Gamma-Y, x),
$$
where $\mathcal{Y}_{u}$ is the set of all signed cycles passing through $u$, and $\Gamma-Y$ is the graph obtained from $\Gamma$ by deleting $Y$.
\end{lemma}
A signed unicyclic graph has the pairing property if and only if its underlying graph is bipartite because it has a unique cycle. Next, we characterize the unique unbalanced bipartite unicyclic graphs with the maximum Estrada index among all unbalanced bipartite unicyclic graphs.
\begin{figure}
\centering
	\includegraphics[scale=.8]{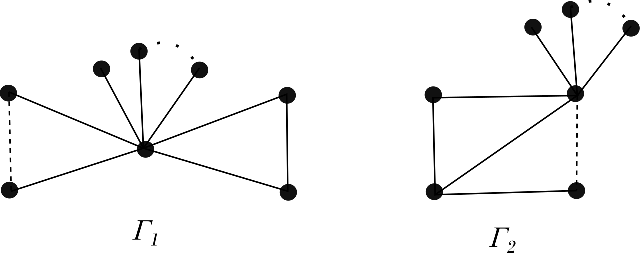}
	\caption{Signed graphs $\Gamma_1$ and  $\Gamma_2$, which are in the statement of Theorem \ref{3.7}.}
	\label{Figure 1}
\end{figure}

\begin{theorem}\label{3.6}
Let $\Gamma_{p}^-(n, n)$ be  the set of all unbalanced bipartite unicyclic graphs on $n$ vertices. If $\Gamma \in \Gamma_{p}^-(n, n)$, then  $EE(\Gamma) \leq EE(\Gamma_n^{4-})$ with equality
if and only if $\Gamma$ is switching equivalent to $ \Gamma_n^{4-}$.
\end{theorem}
{\bf Proof.} Let $\Gamma \in \Gamma_{p}^-(n, n)$ be an unbalanced bipartite unicyclic graph on $n$ vertices. By applying Lemma \ref{3.3}, we get $\mu_1(\Gamma) \leq \mu_1(\Gamma_n^{4-})$ with equality
if and only if $\Gamma$ is switching equivalent to $ \Gamma_n^{4-}$. Now, by Lemma \ref{3.5}, the characteristic polynomial of   $\Gamma_n^{4-}$  is given by
\begin{center}
$\varphi(\Gamma_n^{4-}, x)=x^{n-4}\{x^4-nx^2+2(n-2)\}.$
\end{center}
Clearly, the signed graph $\Gamma_n^{4-}$ has four non-zero eigenvalues. Let the signed graph $\Gamma$, where $\Gamma \in \Gamma_{p}^-(n, n)$, contains a cycle of length $l\geq 4$. Therefore, the unbalanced cycle $C_{l-}$ is an induced signed subgraph of $\Gamma$. The eigenvalues of  $C_{l-}$ are given by $2 \cos \frac{ (2r+1) \pi}{l}, \quad r=0, 1,2, \ldots, l-1$. Since $l$ is a positive even integer and thus all the eigenvalues of  $C_{l-}$ are non-zero. Hence the result follows by Theorem \ref{3.2} and Lemma \ref{3.4}.\qed

\begin{theorem}\label{3.7}
Let $\Gamma_{p}^-(n, n+1)$ be  the set of all unbalanced bicyclic graphs on $n$ $(n\geq5)$ vertices with the pairing property. If $\Gamma \in \Gamma_{p}^-(n, n+1)\backslash \{\Gamma_1, \Gamma_2\}$, $\Gamma$ is not switching equivalent to $ \Gamma_1$ and $ \Gamma_2$,  then $EE(\Gamma_1)> EE(\Gamma_2)>EE(\Gamma), $ where $\Gamma_1$ and $\Gamma_2$ are the signed graphs on $n$ vertices as shown in Fig.2.
\end{theorem}
{\bf Proof.} By Lemma \ref{3.5}, the characteristic polynomials of  $\Gamma_1$ and $\Gamma_2$ are, respectively, given by
\begin{center}
$\varphi(\Gamma_1, x)=x^{n-6}(x^2-1)\{x^4-nx^2+n-5\}$
\end{center}
and
\begin{center}
$\varphi(\Gamma_2, x)=x^{n-4}\{x^4-(n+1)x^2+2(n-2)\}.$
\end{center}
It is easy to see that $Spec(\Gamma_1)=  \{\pm\sqrt{\frac{n\pm\sqrt{n^2-4n+20}}{2}},1,0^{n-6},-1\}$ and $Spec(\Gamma_2)=$ \\ $\{\pm\sqrt{\frac{n+1\pm\sqrt{(n+1)^2-8(n-2)}}{2}},0^{n-4}\}$ respectively. Therefore\\
\begin{equation}
EE(\Gamma_1)=n-6+e^{\sqrt{\frac{n+\sqrt{n^2-4n+20}}{2}}}+e^{\sqrt{\frac{n-\sqrt{n^2-4n+20}}{2}}}+e^{-\sqrt{\frac{n+\sqrt{n^2-4n+20}}{2}}}+e^{-\sqrt{\frac{n-\sqrt{n^2-4n+20}}{2}}}+e+e^{-1}.
\end{equation}
Also,\\
\begin{equation}
\begin{split}
EE(\Gamma_2)&=n-4+e^{\sqrt{\frac{n+1+\sqrt{(n+1)^2-8(n-2)}}{2}}}+e^{\sqrt{\frac{n+1-\sqrt{(n+1)^2-8(n-2)}}{2}}}+e^{-\sqrt{\frac{n+1+\sqrt{(n+1)^2-8(n-2)}}{2}}}\\&+e^{-\sqrt{\frac{n+1-\sqrt{(n+1)^2-8(n-2)}}{2}}}.
\end{split}
\end{equation}
Note that $\sqrt{\frac{n+\sqrt{n^2-4n+20}}{2}}> \sqrt{\frac{n+1+\sqrt{(n+1)^2-8(n-2)}}{2}}$ for $n \geq 5$. We can check that the right-hand side of $(3.18)$ is greater than of $(3.19)$ for $n\geq 5$, which proves that $EE(\Gamma_1)> EE(\Gamma_2)$.\\
\indent Let  $\Gamma \in \Gamma_{p}^-(n, n+1)\backslash \{\Gamma_1, \Gamma_2\}$ be an unbalanced bicyclic graph  with the pairing property. To prove the result it is enough to show that $EE(\Gamma_2)>EE(\Gamma)$, where $\Gamma$ is not switching equivalent to  $ \Gamma_1$ and $ \Gamma_2$. Since $\mu_1(\Gamma_1)> \mu_1(\Gamma_2)$, for $n\geq 5$, therefore, by [\cite{bo}, Theorem 3.1], we get $\mu_1(\Gamma_1)>\mu_1(\Gamma_2)> \mu_1(\Gamma)$, where $\Gamma$ is not switching equivalent to $ \Gamma_1$ and $ \Gamma_2$. The signed graph $\Gamma_2$ has four non-zero eigenvalues. Clearly, the signed graph $\Gamma$ contains a 4-vertex signed path $P_4$ as an induced subgraph for $n\geq 5$. The characteristic polynomial of $P_4$ is given by
\begin{center}
$\varphi(P_4, x)=x^4-3x^2+1.$
\end{center}
Thus the signed path $P_4$ has four non-zero eigenvalues. By interlacing Lemma \ref{3.4},  the signed graph $\Gamma$ has at least four non-zero eigenvalues. Hence the result follows by Theorem \ref{3.2}. \qed

\begin{lemma} \label{3.8}{\em \cite{z1}}
Let $(G, \sigma)$ be a connected signed graph. Then, we have $\mu_{1}((G, \sigma)) \leq \mu_{1}((G, +))$ with equality if and only if $(G, \sigma)$ switches to $(G, +)$.
 \end{lemma}

\begin{lemma} \label{3.9} {\em \cite{z2}}
 For an eigenvalue $\mu$ of a connected signed graph $\Gamma$, there exists a switching equivalent signed graph $\Gamma'$, for which the $\mu$-eigenspace contains an eigenvector whose non-zero coordinates are of the same sign.
 \end{lemma}
  Let $\Gamma = (K_{m, n},\sigma)$ be a  signed complete bipartite graph, where  $K_{m, n}$ is  the complete bipartite graph on $m+n$ vertices. We denote, $S(K_{m, n},-)$, the set of all unbalanced complete bipartite graphs on $n+m$ vertices. Also, let $\Gamma^*$ be an unbalanced complete bipartite graph that contains exactly one negative edge. Finally, we characterize the unbalanced complete bipartite graphs with the maximum Estrada index.
 \begin{theorem}
  Let $S(K_{m, n},-)$ be the set of all unbalanced complete bipartite graphs on $n+m$ vertices. If $\Gamma \in S(K_{m, n},-)$, then $EE(\Gamma^*)\geq EE(\Gamma) $ with equality if and only if $\Gamma$ is switching equivalent to $\Gamma^*$, where $\Gamma^*$ is an unbalanced complete bipartite signed graph that contains exactly one negative edge.
 \end{theorem}
 {\bf Proof.} Let $\left(\Gamma_{1}, \Gamma_{2}, \ldots, \Gamma_{k}\right)$ be a sequence  consisting of the representatives of all switching equivalence classes of unbalanced  complete bipartite graphs with $m+n$ vertices such that the representatives are ordered non-increasingly by the largest eigenvalue (index) and chosen in such a way that, for $1 \leq j \leq k$, the $\mu_{1}(\Gamma_{j})$-eigenspace  contains an eigenvector whose non-zero coordinates are positive (This existence of $\Gamma_{j}$ is provided by Lemma \ref{3.9}). By Lemma \ref{3.8}, the signed complete bipartite with all positive signature has the maximum index. Therefore, by [\cite{z3},Theorem $3.2$], the  signed graph with maximum index among all unbalanced complete bipartite graphs on $n+m$ vertices is $\Gamma^*$ (upto switching), where $\Gamma^*$ is an unbalanced complete bipartite signed graph that contains exactly one negative edge. That is, for each $\Gamma \in S(K_{m, n},-)$, we have $\mu_1(\Gamma) \leq \mu_1(\Gamma^*)$ with equality if and only if $\Gamma$ is switching equivalent to $\Gamma^*$. Now, by [\cite{ts}, Theorem 4.1], the spectrum of the signed graph $\Gamma^*$ is given by
 $$Spec(\Gamma^*)=\{\pm\sqrt{\frac{mn \pm \sqrt{n^2+2(m-1)(n-2)^2+n^2(m-1)^2}}{2}},0^{n-4}\}.$$ Clearly, the signed graph $\Gamma^*$ has four non-zero eigenvalues. Also, with a suitable labelling of the vertices of $\Gamma \in S(K_{m, n},-)$,  its adjacency matrix is given by  $$ A_\Gamma
=\left(\begin{array}{cc}
O_{m\times m} & B_{m\times n} \\
 B_{n\times m } ^{\top}& O_{n\times n}
\end{array}\right),
$$ where $B_{m\times n}$ is a matrix whose entries are either $1$ or $-1$.  We know that $rank(A_\Gamma)= rank(B_{m\times n} )+rank(B_{m\times n}^\top )=2 rank(B_{m\times n} )$. It is easy to see that $rank(A_\Gamma)\geq 4$ because if $rank(B_{m\times n} )=1$, then  $\Gamma$ is a switching equivalent to a signed complete bipartite graph with all positive signature which is a contradiction. Thus the signed graph $\Gamma \in S(K_{m, n},-)$ has at least four non-zero eigenvalues. Hence the result follows by Theorem \ref{3.2}.\qed

\noindent{\bf Acknowledgements.}   This research is supported by SERB-DST research project number CRG/2020/000109. The research of Tahir Shamsher is supported by SRF financial assistance by Council of Scientific and Industrial Research (CSIR), New Delhi, India.\\

\noindent{\bf Data availibility} Data sharing is not applicable to this article as no datasets were generated or analyzed
during the current study.\\

\noindent{\bf Declaration}
The authors declare that there is no conflict of interest amongst them.

\end{document}